\documentclass[titlepage,twoside,12pt]{article}
\usepackage{amsmath}
\usepackage{amssymb}
\usepackage{amsfonts}
\begin{document}

{\LARGE \bf Dynamics in the Category {\bf Set}} \\ \\

{\bf Elem\'{e}r E ~Rosinger} \\ \\
{\small \it Department of Mathematics \\ and Applied Mathematics} \\
{\small \it University of Pretoria} \\
{\small \it Pretoria} \\
{\small \it 0002 South Africa} \\
{\small \it eerosinger@hotmail.com} \\ \\

\hfill {\it Dedicated to Marie-Louise Nykamp} \\ \\

{\bf Abstract} \\

What makes sets, or more precisely, the category {\bf Set} important in Mathematics are the well known {\it two} specific ways in which arbitrary mappings $f : X \longrightarrow Y$ between any two sets $X, Y$ can {\it fail} to be bijections. Namely, they can fail to be injective, and/or to be surjective. As for bijective mappings they are rather trivial, since with some relabeling of their domains or ranges, they simply become permutations, or even identity mappings. \\
To the above, one may add the {\it third} property of sets, namely that, between any two nonvoid sets there exist mappings. \\
These three properties turn out to be at the root of much of the interest which the category {\bf Set} has in Mathematics. Specifically, these properties create a certain {\it dynamics}, or for that matter, lack of it, on the level of the category {\bf Set} and of some of its subcategories. \\ \\ \\

{\bf 0. Preliminaries} \\

Let us consider in the category {\bf Set} the subcategory ${\bf Set_{\bf bij}}$ where the objects are all the sets $X$ in {\bf Set}, while the morphisms are only the bijections $f : X \longrightarrow Y$ between sets $X, Y$ in {\bf Set}. Further, let us consider the subcategory ${\bf Set_{\bf id}}$ in ${\bf Set_{\bf bij}}$ where the objects are again all the sets $X$ in {\bf Set}, while the morphisms are only the identity mappings $id_X : X \longrightarrow X$. \\

One may at first sight possibly comment that in the increasing chain of categories \\

(0.1)~~~~ $ {\bf Set_{\bf id}} ~\subset~ {\bf Set_{\bf bij}} ~\subset~ {\bf Set} $ \\

in the left term category nothing at all happens, in the middle one the only thing happening is re-labelings and permutations, while just about everything else that happens to happen in the Mathematics based on Set Theory does in fact take place in the right term category. \\
In other words, one may wonder whether not all Mathematics based on Set Theory does actually take place due to the fact that there are mappings in the category {\bf Set} other than mere bijections. \\

Let us therefore address the most simple basics of the difference within the category {\bf Set} between bijections, and on the other hand, all the other mappings which are not bijections. \\

Given any mapping \\

(0.2)~~~~ $ f : X \longrightarrow Y $ \\

in ${\bf Set}$, one defines the equivalence relation $\approx_f$ on $X$ by $x \approx_f x\,' \Leftrightarrow f ( x ) = f ( x\,' )$. Further, one defines $X_f = X / \approx_f$, and lastly, one defines $\underline f : X_f \longrightarrow f ( X )$ by $\underline f ( ( x )_{\approx_f}) = f ( x )$, for $x \in X$, where $( x )_{\approx_f}$ denotes the $\approx_f$ equivalence class of $x$. Then obviously \\

(0.3)~~~~ $ X \stackrel{i_f}\longrightarrow X_f \stackrel{\underline f}\longrightarrow f ( X )
                                                         \stackrel{i}\longrightarrow Y $ \\

where $i_f$ is surjective, $\underline f$ is bijective, and $i$ is injective. \\

In this way one obtains a mapping $\rho : {\bf Set} \longrightarrow {\bf Set_{\bf bij}}$ given by \\

(0.4)~~~~ $ ( f : X \longrightarrow Y ) \stackrel{\rho}\longrightarrow
                                 ( \underline f : X_f \longrightarrow f ( X ) ) $ \\

{\bf Remark 0.1.} \\

The mapping $\rho$ in (0.4) is {\it not} a functor, since both the domain $X_f$ and range $f ( X )$ of the mapping $\underline f$ may depend on the mapping $f$, and not only on $X$ and $Y$, respectively. \\

{\bf Problem 0.1.} \\

Based on (0.4), study the difference between a {\it bijection}, and on the other hand, an {\it arbitrary mapping}. Namely, do so by studying the difference between $f : X \longrightarrow Y$ and $\underline f : X_f \longrightarrow f ( X )$, then classify it, and find conditions in terms of such classification when such a difference does not exists, that is, when $f = \underline f$.

\hfill $\Box$ \\

A natural way to proceed is to make the following two comparisons \\

(0.5)~~~~ between $ X $ and $ X_f $ \\

and then also \\

(0.6)~~~~ between $ \phi $ and $ Y \setminus f ( X ) $ \\

Let us therefore associate with (0.2) the following entity \\

(0.7)~~~~ $ dev ( f ) = ( X_f, \, Y \setminus f ( X ) ) $ \\

called the {\it deviation} of $f$. \\

Clearly \\

(0.8)~~~~ $ dev ( f ) =  ( X, \phi ) ~\Longleftrightarrow~ f $ bijective \\

(0.9)~~~~ $ dev ( f ) =  ( X_f, \phi ) ~\Longleftrightarrow~ f $ surjective \\

(0.10)~~~~ $ dev ( f ) =  ( X, Y \setminus f ( X ) ) ~\Longleftrightarrow~ f $ injective \\

where we make the following identifications \\

(0.11)~~~~ $ \begin{array}{l}
                      X = X_{id_X} = \{ \{ x \} ~|~ x \in X \} \\ \\
                      X_f = X / \approx_f ~=~ \{ f^{-1} ( \{ y \} ) ~|~ y \in f ( X ) \}
             \end{array} $ \\ \\

{\bf Remark 0.2.} \\

1) In view of (0.1), one may consider that, given a subcategory ${\bf C}$, with \\

(0.12)~~~~ $ {\bf Set_{\bf bij}} ~\subset~ {\bf C} ~\subset~ {\bf Set} $ \\

then the larger the deviations $dev ( f )$ of mappings $f$ in ${\bf C}$ from bijections, the more dynamism in the subcategory ${\bf C}$. \\

2) Obviously, when measuring the deviation of an arbitrary mapping $f : X \longrightarrow Y$ in the category ${\bf Set}$, we cannot use as term of reference identity mappings, since they only cover the particular cases when the domains of the mappings are the same with their ranges, that is, $X = Y$. Consequently, deviations of arbitrary mappings in the category ${\bf Set}$ should rather refer to bijections. \\ \\

{\bf 1. Comparing Deviations} \\

 A natural way to compare the second terms $Y \setminus f ( X )$ in the deviations (0.7) is obvious, namely, through the usual inclusion relation $\subseteq$ between subsets of $Y$. As for the comparison of the first terms $X_f$, this can equally naturally be done as follows. We consider the set ${\cal PR} ( X )$ of all {\it partitions} of $X$. And then in view of (0.11), we obtain \\

(1.1)~~~~ $ X_f \in {\cal PR} ( X ) $ \\

for all mappings $f : X \longrightarrow Y$ in ${\bf Set}$. Further, we consider on ${\cal PR} ( X )$ the partial order $\leq$ defined by \\

(1.2)~~~~ $ {\cal A} \leq {\cal B} ~\Longleftrightarrow~
                         \forall~ A \in {\cal A} ~:~ \exists~ B \in {\cal B} ~:~ A \subseteq B $ \\

It follows in particular that, see (0.8)-(0.10) \\

(1.3)~~~~ $ X_{id_X} \leq X_f \leq \{ X \} $ \\

for all mappings $f : X \longrightarrow Y$ in ${\bf Set}$, with \\

(1.4)~~~~ $ X_{id_X} = X_f ~\Longleftrightarrow~ f $ injective \\

(1.5)~~~~ $ X_f = \{ X \} ~\Longleftrightarrow~ f $ constant \\

And now, given two mappings $f : X \longrightarrow Y$ and $f\,' : X \longrightarrow Y$ in ${\bf Set}$, we can compare their deviations according to the following definition \\

(1.6)~~~~ $ dev ( f ) \leq dev( f\,' ) ~\Longleftrightarrow~
                        X_f \leq X_{f\,'},~~ Y \setminus f ( X ) \subseteq Y \setminus f\,' ( X ) $ \\

In view of (0.7), we obtain \\

{\bf Lemma 1.1.} \\

Given a mapping $f : X \longrightarrow Y$ in ${\bf Set}$, then \\

(1.7)~~~~ $ f $ is bijective $ ~\Longleftrightarrow~ dev ( f ) = ( X, \phi ) \leq dev ( f\,' ) $ \\

for all mappings $f\,' : X \longrightarrow Y$ in ${\bf Set}$.

\hfill $\Box$ \\

Let us now see the way deviations behave with respect to the composition of mappings. It will be convenient to consider separately the two components in the deviations defined in (0.7), namely, to denote \\

(1.8)~~~~ $ dev_1 ( f ) = X_f,~~ dev_2 ( f ) = Y \setminus f ( X ) $ \\

for any mapping $f : X \longrightarrow Y$ in ${\bf Set}$. \\

{\bf Theorem 1.1.} \\

Let $X \stackrel{f}\longrightarrow Y \stackrel{g}\longrightarrow Z$ in ${\bf Set}$, and $X \stackrel{h}\longrightarrow Z$, with $h = g \circ f$. Then \\

(1.9)~~~~ $ dev_1 ( f ) \leq dev_1 ( h ) $ \\

{\bf Proof.} \\

According to (1.8), we have \\

(1.10)~~~~ $ dev_1 ( f ) = X_f,~~ dev_2 ( h ) = X_h $ \\

while (0.11) gives \\

(1.11)~~~~ $ X_f = \{ f^{-1} ( \{ y \} ) ~|~ y \in f ( X ) \},~~ X_h = \{ h^{-1} ( \{ y \} ) ~|~ y \in h ( X ) \} $ \\

Take $A \in X_f$, then $A = f^{-1} \{ y \}$, for some $y \in Y$. Take $z = g ( y )$, then \\

$~~~~ h^{-1} ( \{ z \} ) = ( g \circ f)^{-1} ( \{ z \} ) =
                         f^{-1} ( g^{-1} ( \{ z \} ) ) \supseteq f^{-1} ( \{ y \} ) = A $

\hfill $\Box$ \\

As for the behaviour with respect to composition of mappings of the second component of deviations, the situation is more complicated, namely \\

{\bf Theorem 1.2.} \\

Let $X \stackrel{f}\longrightarrow Y \stackrel{g}\longrightarrow Z$ in ${\bf Set}$, and $X \stackrel{h}\longrightarrow Z$,
with $h = g \circ f$. Then the only general relationship regarding $dev_2$ is \\

(1.12)~~~~ $ dev_2 ( g ) \leq dev_2 ( h ) $ \\

{\bf Proof.} \\

In view of (1.10), we have \\

$~~~~ dev_2 ( g ) = Z \setminus g ( Y ),~~ dev_2 ( h ) = Z \setminus h ( X ) $ \\

and obviously $h ( X ) \subseteq g ( Y )$. \\

To complete the proof, we give examples that $dev_2 ( f )$ cannot in general be compared in any given fixed way with either terms in (1.12), even in the particular case when $X = Y = Z$. \\

Let $f$ be a surjection, then $dev_2 ( f ) = \phi$. Thus the only fixed way it could in general relate to the two terms in (1.12) would be by \\

(1.13)~~~~ $ dev_2 ( f ) = \phi \leq dev_2 ( g ),~~ dev_2 ( f ) = \phi \leq dev_2 ( h ) $ \\

On the other hand, if $g$ is a surjection, then the only fixed way it could in general relate to $dev_2 ( f )$ would be \\

(1.14)~~~~ $ dev_2 ( g ) = \phi \leq dev_2 ( f ) $ \\

and the last two relations contradict one another in case of arbitrary mappings $f$ and $g$, even when $X = Y = Z$. \\

{\bf Remark 1.1.} \\

The results in Theorems 1.1. and 1.2. above indicate the {\it complexity} in the behaviour of deviations under the composition of functions. In other words, and according to customary intuition, the deviation of mappings from being bijective increases through their composition. \\
And it is precisely in this complexity in the behaviour of deviations from bijectivity under the composition of functions that one can see the source of much of the Mathematics based on Set Theory. \\ \\

{\bf 2. Deviations in Abelian Groups} \\

As seen in Theorem 1.2., the second component $dev_2$ in deviations has a complex behaviour with respect to the composition of mappings in the category ${\bf Set}$. Consequently, it may be useful to consider certain alternative definitions of it in suitable subcategories of the category ${\bf Set}$. \\
Here we do that in the category ${\bf AG}$ of Abelian groups. \\

Let $f : X \longrightarrow Y$ be a group homomorphism in the category ${\bf AG}$, then we define its {\it group deviation} by \\

(2.1)~~~~ $ devg ( f ) = ( X_f, Y / f ( X ) ) $ \\

as well as its respective components \\

(2.2)~~~~ $ devg_1 ( f ) = dev ( f ) = X_f = X / ker ( f ),~~ devg_2 ( f ) = Y / f ( X ) $ \\

Obviously, both $devg_1 ( f )$ and $devg_2 ( f )$ are Abelian groups. Furthermore, as a corresponding modification of the partial order (1.2), we consider the following one. Given two group homomorphisms $f : X \longrightarrow Y$ and $f\,' : X \longrightarrow Y$ in ${\bf AG}$, then  \\

(2.3)~~~~ $ devg ( f ) \leq devg ( f\,' ) ~\Longleftrightarrow~
               \left ( ~ \begin{array}{l}
                                devg_1 ( f ) ~~\mbox{subgroup of}~~ devg_1 ( f\,' ) \\ \\
                                devg_2 ( f ) ~~\mbox{subgroup of}~~ devg_2 ( f\,' )
                          \end{array} ~ \right ) $ \\

{\bf Lemma 2.1.} \\

Given a group homomorphism $f : X \longrightarrow Y$ in ${\bf AG}$, then \\

(2.4)~~~~ $ f $ group isomorphims $ ~\Longleftrightarrow~ devg ( f ) = ( X, 0 ) \leq devg ( f\,' ) $ \\

for all group homomorphisms $f\,' : X \longrightarrow Y$ in ${\bf AG}$, where $0$ denotes the trivial subgroup in the Abelian group $Y$. Also \\

(2.5)~~~~ $ f $ surjective $ ~\Longleftrightarrow~ devg ( f ) = ( X_f, 0 ) $ \\

(2.6)~~~~ $ f $ injective $ ~\Longleftrightarrow~ devg ( f ) = ( X, Y / f ( X ) ) $

\hfill $\Box$ \\

And now, let us see the way the group deviations (2.1) behave with respect to the composition of group
homomorphisms. \\

{\bf Theorem 2.1.} \\

Let $X \stackrel{f}\longrightarrow Y \stackrel{g}\longrightarrow Z$ group homomorphisms in ${\bf AG}$, and $X \stackrel{h}\longrightarrow Z$, with $h = g \circ f$. Then \\

(2.7)~~~~ $ devg_1 ( f ) \leq devg_1 ( h ) $ \\

And the only general relationship regarding $devg_2$ is \\

(2.8)~~~~ $ devg_2 ( g ) $ is a subgroup of $ devg_2 ( h )$ \\

{\bf Proof.} \\

The relation (2.7) follows from (2.2) and Theorem 1.1.

As for (2.8), let us assume the particular case when we have the identity $X = Y = Z$ of the three Abelian groups. Further, let $f$ be a surjective group homomorphism, then (2.5) gives $devg ( f ) = ( X_f, 0 )$, thus $devg_2 ( f ) = 0$. It follows that the only fixed way $devg_2 ( f )$ could in general relate to the two terms in (2.8) would be by \\

(2.9)~~~~ $ \begin{array}{l}
                dev_2 ( f ) = 0 $ is a subgroup of $ dev_2 ( g ) \\ \\
                dev_2 ( f ) = 0 $ is a subgroup of $ dev_2 ( h )
             \end{array} $ \\

On the other hand, if $g$ is a surjective group homomorphism, then the only fixed way $devg_2 ( g ) = 0$ could in general relate to $dev_2 ( f )$ would be \\

(2.10)~~~~ $ dev_2 ( g ) = 0 $ is a subgroup of $ dev_2 ( f ) $ \\

and the last two relations contradict one another in case of arbitrary group homomorphisms $f$ and $g$, even when we have three identical Abelian groups $X = Y = Z$. \\

{\bf Remark 2.1.} \\

In view of Theorem 2.1. above, the transition from the category ${\bf Set}$ to its subcategory ${\bf AB}$, as well as from the deviations (0.7) to the group deviations (2.1) does not change the fact that $devg_2$ has again a complex behaviour, similar with $dev_2$, this time with respect to the composition of group homomorphisms. Needless to say, this fact is not surprising when one thinks that Group Theory may have much of its source precisely in this complex behaviour of $devg_2$, just as the Mathematics base on Set Theory seems to have most of its source in the complex behaviour of $dev_2$, see Remark 1.1. \\ \\

{\bf 3. The Case of Chu Spaces} \\

Recently, Chu spaces have known a wider interest in Mathematics and Physics, see [1] and the references cited there. \\
For convenience, here is a brief presentation of the concept of Chu space, more precisely, of the various categories of
such spaces. \\

Being given a fixed set $W$, the corresponding category ${\bf Chu}_{\bf W}$ of Chu spaces is defined as follows. The objects of that category are of the form \\

(3.1)~~~~ $ ( X, Y, f ) $ \\

where $X, Y$ are sets, while $f : X \times Y \longrightarrow W$. \\

The morphisms of the category ${\bf Chu}_{\bf W}$ are of the form \\

(3.2)~~~~ $ m = ( \underline m, \overline m ) : ( X, Y, f ) \longrightarrow ( U, V, g ) $ \\

where \\

(3.3)~~~~ $ \underline m : X \longrightarrow U,~~ \overline m : V \longrightarrow Y $ \\

and they are required to satisfy the condition \\

(3.4)~~~~ $ f ( x, \overline m ( v ) ) = g ( \underline m ( x ), v ),~~ x \in X, v \in V $ \\

Lastly, the composition of two morphisms ~$m : ( X, Y, f ) \longrightarrow ( U, V, g )$ and $n : ( U, V, G ) \longrightarrow ( P, Q, h )$ is given by \\

(3.5)~~~~ $ k = m \circ n =
               (~ \underline k = \underline m \circ \underline n,~ \overline k = \overline m \circ \overline n ~) $ \\

We note that for $W = \phi$, or $W$ with one single element, the Chu spaces in ${\bf Chu}_{\bf W}$ are trivial. Therefore,
we shall assume that $\{ 0, 1 \} \subseteq W$. \\

A remarkable property of Chu spaces is that the category ${\bf Set}$ can be {\it fully embedded} into them by the
functor $E : {\bf Set} \longrightarrow {\bf Chu_W}$ defined as follows, [1] \\

(3.6)~~~~ $ ( f : X \longrightarrow Y ) \longmapsto
                 ( ( f, f^{-1} ) : ( X, {\cal P} ( X ), e_X ) \longrightarrow ( Y, {\cal P} ( Y ), e_Y ) ) $ \\

where, for $x \in X, A \subseteq X$, we have \\

(3.7)~~~~ $ e_X ( x, A ) ~=~ \begin{array}{|l}
                                     ~ 1 ~\mbox{if}~ x \in A \\ \\
                                     ~ 0 ~\mbox{if}~ x \notin A
                             \end{array} $ \\

and similarly for $e_Y$. \\

In view of this representation of the category ${\bf Set}$ in Chu spaces, we shall consider the deviation introduced in (0.7) in the alternative terms of Chu spaces. \\

In this regard, first, let us consider the deviation of the above mapping $e_X$ which obviously does not depend on $f$. A similar situation applies then to $e_Y$. We obviously have \\

(3.8)~~~~ $ e_X : ( X \times {\cal P} ( X ) ) \ni ( x, A ) \longmapsto e_X ( x, A ) \in W $ \\

and thus \\

(3.9)~~~~ $ dev_1 ( e_X ) = ( X \times {\cal P} ( X ) )_{e_X},~~~
                      dev_2 ( e_X ) = W \setminus e_X ( X \times {\cal P} ( X ) )$ \\

Without loss of generality, we can assume \\

(3.10)~~~~ $ W = \{ 0, 1 \},~~~ car ( X ) \geq 2 $ \\

Then obviously \\

(3.11)~~~~ $ dev_2 ( e_X ) = \phi $ \\

thus all the information in the deviation $dev ( e_X )$ is contained in $dev_1 ( e_X ) = ( X \times {\cal P} ( X ))_{e_X}
= ( X \times {\cal P} ( X )) / \approx_{e_X}$. In this regard we note that, given $( x, A ), ( x\,', A\,' ) \in X \times {\cal P}
( X )$, we have \\

(3.12)~~~~ $ ( x, A ) \approx_{e_X} ( x\,', A\,' ) ~~~\Longleftrightarrow~
                    \left ( \begin{array}{l}
                                  ~\mbox{either}~ x \in A ~\mbox{and}~ x\,' \in A\,' \\ \\
                                  ~\mbox{or}~ x \notin A ~\mbox{and}~ x\,' \notin A\,'
                             \end{array} \right ) $ \\

and in view of (0.11), we have \\

(3.13)~~~~ $ dev_1 ( e_X ) = \{ {\cal X}_0, {\cal X}_1 \} $ \\

where \\

(3.14)~~~~ $ \begin{array}{l}
                     {\cal X}_0 = \{ ( x, A ) \in X \times {\cal P} ( X ) ~|~ x \notin A \} \\ \\
                     {\cal X}_1 = \{ ( x, A ) \in X \times {\cal P} ( X ) ~|~ x \in A \}
             \end{array} $ \\

Now we return to (3.6) and consider the possible additional information which the respective representation of the
category ${\bf Set}$ in the category ${\bf Chu_{\{ 0, 1 \}}}$ can give on the deviation (0.7) of mappings in ${\bf Set}$. \\
The possibility that (3.6) may provide such further information comes from the fact that, in (3.6), not only the mappings $f$ are involved, but also their inverses $f^{-1}$ which, as seen next, are in fact inverses of certain {\it extensions} $\widetilde f$ of the mappings $f$. \\
Indeed, given in the category ${\bf Set}$ any mapping \\

(3.15)~~~~ $ f : X \longrightarrow Y $ \\

its inverse is the mapping \\

(3.16)~~~~ $ f^{-1} : {\cal P} ( Y ) \longrightarrow {\cal P} ( X ) : U \longmapsto f^{-1} ( U ) $ \\

thus in fact $f^{-1}$ is rather nearer to the inverse of the following extension of $f$ \\

(3.17)~~~~ $ \widetilde f  : {\cal P} ( X ) \longrightarrow {\cal P} ( Y ) : A \longmapsto f ( A ) $ \\

For simplicity of notation we shall write $f$ instead of $\widetilde f$, whenever no confusion may arise. The connection between $f, \widetilde f$ and $f^{-1}$ will be further considered in the sequel, and among others, in (3.58) - (3.61) below. Here we note that, for $A \subseteq X$, we have \\

(3.18)~~~~ $ \widetilde f ( A ) = \phi ~\Longleftrightarrow~ A = \phi $ \\

Given now the deviation of $f : X \longrightarrow Y$, namely \\

(3.19)~~~~ $ dev ( f ) = ( X_f, Y \setminus f ( X ) ) $ \\

let us also consider the deviations of $\widetilde f  : {\cal P} ( X ) \longrightarrow {\cal P} ( Y )$ and $f^{-1} : {\cal P} ( Y ) \longrightarrow {\cal P} ( X )$, which are given by \\

(3.20)~~~~ $ dev ( \widetilde f ) =
                    ( ( {\cal P} ( X ) )_{\widetilde f}, {\cal P} ( Y ) \setminus \widetilde f ( {\cal P} ( X ) ) $ \\

(3.21)~~~~ $ dev ( f^{-1} ) = ( ( {\cal P} ( Y ) )_{f^{-1}}, {\cal P} ( X ) \setminus f^{-1} ( {\cal P} ( Y ) ) ) $ \\

{\bf Lemma 3.1.} \\

(3.22)~~~~ $ f $ injective $ ~\Longleftrightarrow~ \widetilde f $ injective \\

(3.23)~~~~ $ f $ surjective $ ~\Longleftrightarrow~ \widetilde f $ surjective \\

(3.24)~~~~ $ f $ bijective $ ~\Longleftrightarrow~ \widetilde f $ bijective \\

{\bf Proof.} \\

The relation (3.22) is obvious. \\

The implication ''$\Longrightarrow$'' in (3.23) is obvious. \\
The converse implication in (3.23). Let $y \in Y$, then there exists $A \subseteq X$, such that $f ( A ) = \{ y \}$. But (3.18) gives  $A \neq \phi$, since $\{ y \} \neq \phi$. Thus $f ( x ) = y$, for some $x \in A$. \\

Now (3.24) follows from (3.22), (3.23). \\

{\bf Lemma 3.2.} \\

(3.25)~~~~ $ f^{-1} $ surjective $~\Longleftrightarrow~ f $ injective \\

(3.26)~~~~ $ f^{-1}|_{{\cal P} ( Y \cap f ( X ) )} $ injective \\

(3.27)~~~~ $ f^{-1} $ injective $~\Longleftrightarrow~ f $ surjective \\

(3.28)~~~~ $ f^{-1} $ bijective $~\Longleftrightarrow~ f $ bijective \\

{\bf Proof.} \\

We note that \\

(3.29)~~~~ $ \begin{array}{l}
                A \subseteq X \Longrightarrow A \subseteq f^{-1} ( f ( A ) ) \\ \\
                U \subseteq Y \cap f ( X ) \Longrightarrow f ( f^{-1} ( U )) = U
            \end{array} $ \\

and for $f$ injective we have \\

(3.30)~~~~ $ A \subseteq X \Longrightarrow A = f^{-1} ( f ( A ) ) $ \\

The implication ''$\Longleftarrow$'' in (3.25). Take $A \subseteq X$, then $U = f ( A ) \subseteq Y$. Now $f$ injective and (3.30) imply $A = f^{-1} ( U )$. \\
The converse implication in (3.25). Assume $x, x\,' \in X$, such that $x \neq x\,', f ( x ) = f ( x\,' )$, then $U = \{ f ( x ) \} \subseteq Y$ and $\{ x \} \subsetneqq \{ x, x\,' \} \subseteq f^{-1} ( U )$. Thus $\{ x \} \notin f^{-1} ( {\cal P} ( Y ) )$, which contradicts the hypothesis. \\

The relation (3.26). Let $U, V \subseteq Y \cap f ( X ), U \neq V, f^{-1} ( U ) = f^{-1} ( V )$. Then (3.29) gives $U = f ( f^{-1} ( U ) ) = f ( f^{-1} ( V) ) = V$, which is absurd. \\

The implication ''$\Longleftarrow$'' in (3.27) follows from (3.26). \\
The converse implication in (3.27). Take $y \in Y \setminus f ( X )$, then $f^{-1} ( f ( X ) ) = f^{-1} ( f ( X ) \cup \{ y \} )$, which is absurd. \\

The relation (3.28) follows from (3.25), (3.27). \\

{\bf Lemma 3.3.} \\

(3.31)~~~~ $ ( {\cal P} ( Y \cap f ( X ) ) )_{f^{-1}|_{{\cal P} ( Y \cap f ( X ) )}} =
                                                               {\cal P} ( Y \cap f ( X ) ) ) $ \\

(3.32)~~~~ $ f $ surjective $~\Longleftrightarrow~ ( {\cal P} ( Y ) )_{f^{-1}} = {\cal P} ( Y ) $ \\

{\bf Proof.} \\

In view of (0.11) \\

(3.33)~~~~ $ ( X )_f = X ~\Longleftrightarrow~ f $ injective \\

When applied to the mapping \\

(3.34)~~~~ $ f^{-1}|_{{\cal P} ( Y \cap f ( X ) )} : {\cal P} ( Y \cap f ( X ) ) \longrightarrow  {\cal P} ( X )$ \\

we obtain (3.31), in view of (3.26), namely \\

(3.35)~~~~ $ ( {\cal P} ( Y \cap f ( X ) ) )_{f^{-1}|_{{\cal P} ( Y \cap f ( X ) )}} = \\ \\
                  \hspace*{3.5cm} = \{~ ( U )_{f^{-1}} ~|~ U \subseteq Y \cap f ( X ) ~\} =
                                                {\cal P} ( Y \cap f ( X ) ) ) $ \\

where for $U \subseteq Y$, we have \\

(3.36)~~~~ $ ( U )_{f^{-1}_*} = \{~ V \subseteq Y \cap f ( X ) ~|~ f^{-1} ( V ) = f^{-1} ( U ) ~\} = \\ \\
                  \hspace*{3.5cm} = ( U \cap f ( X ) )_{f^{-1}_*} $ \\

with the abbreviating notation $f^{-1}_* = f^{-1}|_{{\cal P} ( Y \cap f ( X ) )}$. \\

Now (3.32) follows from (3.27). \\

{\bf Lemma 3.3.} \\

(3.37)~~~~ $ f $ injective $~\Longleftrightarrow~ dev ( f^{-1} ) = ( {\cal P} ( Y \cap f ( X ) ), \phi ) $ \\

(3.38)~~~~ $ f $ surjective $~\Longleftrightarrow~ dev ( f^{-1} ) =
                                     ( {\cal P} ( Y ), {\cal P} ( X ) \setminus f^{-1} ( {\cal P} ( Y ) ) $ \\

(3.39)~~~~ $ f $ bijective $~\Longleftrightarrow~ dev ( f^{-1} ) = ( {\cal P} ( Y ), \phi ) $ \\

{\bf Proof.} \\

The implication ''$\Longrightarrow$'' in (3.37). If $f$ is injective, then (3.25), (3.21), (3.16) and (3.31) yield the right term of (3.37). \\
Conversely, $dev_2 ( f^{-1} ) = \phi$ and (3.25) imply that $f$ is injective. \\

The implication ''$\Longrightarrow$'' in (3.38). If $f$ is surjective, then (3.32), (3.21) give (3.38). \\
Conversely, $dev_1 ( f^{-1} ) = {\cal P} ( Y )$ and (3.32) imply that $f$ is surjective. \\

The relation (3.39) is implied by (3.37), (3.38).

\hfill $\Box$ \\

In view of the above, we obtain the following three theorems : \\

{\bf Theorem 3.1.} \\

The following are equivalent \\

(3.40)~~~~ $ f $ injective \\

(3.41)~~~~ $ \widetilde f $ injective \\

(3.42)~~~~ $ f^{-1} $ surjective \\

(3.43)~~~~ $ dev ( f ) = ( X, Y  \setminus f ( X ) ) $ \\

(3.44)~~~~ $ dev ( \widetilde f ) = ( {\cal P} ( X ), {\cal P} ( Y  \setminus f ( X ) ) ) $ \\

(3.45)~~~~ $ dev ( f^{-1} ) = ( {\cal P} ( Y \cap f ( X ) ), \phi ) $ \\

{\bf Theorem 3.2.} \\

The following are equivalent \\

(3.46)~~~~ $ f $ surjective \\

(3.47)~~~~ $ \widetilde f $ surjective \\

(3.48)~~~~ $ f^{-1} $ injective \\

(3.49)~~~~ $ dev ( f ) = ( ( X )_f, \phi ) $ \\

(3.50)~~~~ $ dev ( \widetilde f ) = ( ( {\cal P} ( X ) )_{\widetilde f}, \phi ) $ \\

(3.51)~~~~ $ dev ( f^{-1} ) = ( {\cal P} ( Y ), {\cal P} ( X ) \setminus f^{-1} ( {\cal P} ( Y ) ) $ \\

{\bf Theorem 3.3.} \\

The following are equivalent \\

(3.52)~~~~ $ f $ bijective \\

(3.53)~~~~ $ \widetilde f $ bijective \\

(3.54)~~~~ $ f^{-1} $ bijective \\

(3.55)~~~~ $ dev ( f ) = ( X, \phi ) $ \\

(3.56)~~~~ $ dev ( \widetilde f ) = ( {\cal P} ( X ), \phi ) $ \\

(3.57)~~~~ $ dev ( f^{-1} ) = ( {\cal P} ( Y ), \phi ) $

\hfill $\Box$ \\

Let us note now that, in view of (3.15) - (3.17), (3.29), (3.30) we have \\

(3.58)~~~~ $ f $ injective $ ~\Longrightarrow~ f^{-1} \circ \widetilde f = id_{{\cal P} ( X )} $ \\

(3.59)~~~~ $ \widetilde f \circ f^{-1} = id_{{\cal P} ( Y \cap f ( X ) )} $ \\

(3.60)~~~~ $ f $ surjective $ ~\Longrightarrow~ \widetilde f \circ f^{-1} = id_{{\cal P} ( Y )} $ \\

and we can recall (3.26), namely \\

(3.61)~~~~ $ f^{-1}|_{{\cal P} ( Y \cap f ( X ) )} $ injective \\ \\

{\bf 4. Asymmetries in the Category ${\bf Set}$} \\

The results above show an {\it asymmetry} in the category ${\bf Set}$ when it comes to certain features of arbitrary mappings $f : X \longrightarrow Y$ in that category. And in particular, some of such asymmetries are shown by the deviation $dev$ of such mappings. Here we shall detail two such asymmetries. \\

Clearly, in the categories, see (0.1) \\

(4.1)~~~~ $ {\bf Set_{\bf id}} ~\subset~ {\bf Set_{\bf bij}} $ \\

there are no asymmetries related to the mappings involved. As for the category ${\bf Set}$, we can mention the following asymmetries which can be noted above, and which originate in the {\it two failures} of arbitrary mappings in that category, namely, the failure to be injective, or to be surjective, thus in sort, the failure to be bijective. \\

The {\it injectivity} of a mapping $f : X \longrightarrow Y$ is equivalent with the condition \\

(4.2)~~~~ $ X_f = X_{id_X} $ \\

which occurs within ${\cal P} ( {\cal P} ( X ) )$, while on the other hand, the {\it surjectivity} of that mapping is equivalent with the condition \\

(4.3)~~~~ $ f ( X ) = Y $ \\

which occurs within ${\cal P} ( Y )$, thus one level {\it lower} in the hierarchy \\

(4.4)~~~~ $set,~~ {\cal P} ( set ),~~ {\cal P} ( {\cal P} ( set ) ),~~
                                {\cal P} ( {\cal P} ( {\cal P} ( set ) ) ), \ldots $ \\

A further asymmetry related to mappings $f : X \longrightarrow Y$ in the category ${\bf Set}$ happens with respect to the two components  $dev_1 ( f ) = X_f$ and $dev_2 ( f ) = Y \setminus f ( X )$ of their deviations $dev ( f )$, see (1.8). Namely, as seen in Theorems 1.1. and 1.2., these two components behave differently with respect to the composition of mappings, with the second component having a more complex behaviour. \\

\end{document}